\documentclass[12pt]{article}
\usepackage{cite,amsfonts,amsmath,subeqn}
\def\be{\begin{equation}}
\def\eqn#1{\be\label{#1}}

\def\bea{\begin{eqnarray}}
\def\eqnn#1{\bea\label{#1}}
\def\eea{\end{eqnarray}}

\newcommand{\eqna}[1]{\begin{subequations} \label{#1}
\begin{eqnarray}}
\def\eena{\end{eqnarray}
\end{subequations}}
\def\hR{\hat R} 
  \def\k{\kappa} 
\def\mt{\mapsto}
\def\nn{\nonumber}

\def\cP{{\cal P}} \def\cQ{{\cal Q}}

\newcommand{\RR}{\mbox{${\mathbb R}$}}

\normalbaselines 
\begin{document}

 \pagestyle{empty}

 \begin{center}
 
 \textsf{\LARGE Spectral Decomposition and Baxterisation of  Exotic\\[1mm] 
Bialgebras and Associated Noncommutative\\[2mm] Geometries} 
 
 \vspace{10mm}
 
 {\large D.~Arnaudon$^{a,}$\footnote{Daniel.Arnaudon@lapp.in2p3.fr}, 
 ~A.~Chakrabarti$^{b,}$\footnote{chakra@cpht.polytechnique.fr},\\[2mm] 
 V.K.~Dobrev$^{c,}$\footnote{dobrev@inrne.bas.bg} 
 ~and~ S.G.~Mihov$^{c,}$\footnote{smikhov@inrne.bas.bg}}
 
 \vspace{5mm}
 
 \emph{$^a$ Laboratoire d'Annecy-le-Vieux de Physique Th{\'e}orique LAPTH}
 \\
 \emph{CNRS, UMR 5108, associ{\'e}e {\`a} l'Universit{\'e} de Savoie}
 \\
 \emph{LAPTH, BP 110, F-74941 Annecy-le-Vieux Cedex, France}
 \\
 \vspace{3mm}
 \emph{$^b$ Centre de Physique Th{\'e}orique, CNRS UMR 7644}
 \\
 \emph{Ecole Polytechnique, 91128 Palaiseau Cedex, France.}
 \\
 \vspace{3mm}
 \emph{$^c$ Institute of Nuclear Research and Nuclear Energy} 
 \\
 \emph{Bulgarian Academy of Sciences}
 \\
 \emph{72 Tsarigradsko Chaussee, 1784 Sofia, Bulgaria}
 \\
 \vspace{3mm}

 \end{center}

 \vspace{.8 cm}
 
\begin{abstract}
We study the geometric aspects of two exotic bialgebras ~$S03$~
and ~$S14$~ introduced in ~math.QA/0206053. These bialgebras are
obtained by the Faddeev-Reshetikhin-Takhtajan RTT prescription
with non-triangular R-matrices which are denoted $R_{03}$ and
$R_{14}$ in the classification of Hietarinta, and they are not
deformations of either GL(2) or GL(1/1).  We give the spectral
decomposition which involves two, resp., three, projectors. These
projectors are then used to provide the Baxterisation procedure
with one, resp., two, parameters.  Further, the projectors are
used to construct the noncommutative planes together with the
corresponding differentials following the Wess-Zumino
prescription. In all these constructions there appear
non-standard features which are noted. Such features show the
importance of systematic study of all bialgebras of four
generators.
\end{abstract}
 
 \vfill

\rightline{LAPTH-931/02,\ {}INRNE-TH-02-03}
\rightline{math.QA/0209321}
\rightline{September 2002}
 
 \newpage

\pagestyle{plain}
\setcounter{page}{1}

\section{Introduction}
\label{sect:intro}
\setcounter{equation}{0}
              
Until very recently there was no complete list of the matrix bialgebras 
which are unital associative algebras generated by four elements.    
The list, of course, includes the four cases which are deformations 
of classical ones: two two-parameter deformations  
of each of ~$GL(2)$~ and ~$GL(1|1)$, namely, 
the standard ~$GL_{pq}(2)$\ \cite{DMMZ},  
nonstandard (Jordanian) $GL_{gh}(2)$ \cite{Ag}, 
the standard $GL_{pq}(1|1)$ \cite{HiRi,DaWa,BuTo} 
and the hybrid (standard-nonstandard) 
$GL_{qh}(1|1)$ \cite{FHR}. (Later,  
in \cite{AACDM} it was shown that there are no more 
deformations of $GL(2)$ or $GL(1|1)$.) 
The list includes also five exotic cases which are not deformations of the 
classical algebra of functions over the group $\ GL(2)$ or the 
supergroup $GL(1|1)$. These correspond to $4\times 4$ $R$-matrices 
which are not deformations of the trivial $R$-matrix. 
In the classification of \cite{Hietarinta} there are altogether 
five nonsingular such $R$-matrices. The three triangular ones were 
introduced in \cite{AACDM} and their duals were found and studied in detail in 
\cite{ACDM1}. The study of the two non-triangular cases was started in 
\cite{ACDM2}. There the duals were found and their irreducible representations 
were constructed.  In the present paper we continue the study of the 
non-triangular cases with the  geometric aspects, which are very 
important also for the applications.

\section{$S03$}

\subsection{Spectral decomposition}

We start with the first (of two) nonsingular non-triangular $R$-matrix 
in \cite{Hietarinta} which is not a deformation of the unit matrix:
\begin{equation}
 \label{eq:S03}
 R_{S0,3}\ \equiv\ 
 \left(
 \begin{array}{cccc}
 1 & 0 & 0 & 1 \cr
 0 & 1 & 1 & 0 \cr
 0 & 1 & -1 & 0 \cr
 -1 & 0 & 0 & 1 \cr 
 \end{array}
 \right)
\end{equation} 
The actual tool for the spectral decomposition is the braid matrix 
~$\hat{R} = PR$, where ~$P$~ is the permutation matrix:
\eqn{perm} 
P\ \equiv\ \left(
 \begin{array}{cccc}
 1 & 0 & 0 & 0 \cr
 0 & 0 & 1 & 0 \cr
 0 & 1 & 0 & 0 \cr
 0 & 0 & 0 & 1 \cr       
 \end{array}
 \right)
\end{equation} 
With ~$R=  R_{S0,3}$~ we have for  the braid matrix:
\begin{equation}
  \label{eq:am0}
  \hat{R} = PR =
  \left(
 \begin{array}{cccc}
    1 & 0 & 0 & 1 \cr
    0 & 1 & -1 & 0 \cr
    0 & 1 & 1 & 0 \cr
    -1 & 0 & 0 & 1                                                   
 \end{array}
 \right)     \  .
\end{equation}
We need the minimal polynomial\ pol$(\cdot)$ in
one variable such that \ pol$(\hR)=0$ is the lowest order
polynomial identity satisfied by $\hR$. In the case at hand, 
this identity is:
\begin{equation}
  \label{eq:am1}
 \hat{R}^2 - 2\hat{R} +2I = 0 
\end{equation} 
or
\begin{equation}
  \label{eq:am2}
  (\hat{R} - (1+i)I) (\hat{R} - (1-i)I)=0
\end{equation}

The last identity encodes the projectors we need. Indeed, define:
\begin{equation}
  \label{eq:am3}
  P_{(\pm)} ~\equiv ~ \frac{1}{2} (I \pm i (\hat{R} -I)) 
\end{equation}
Then it is easy to see that ~$P_{(\pm)}$~ satisfy the projector 
properties - orthogonality: 
$$P_{(i)}P_{(j)} = {{\delta}_{ij}}P_{(i)} \ ,\quad i,j = +,- \ ,$$
and resolution of the identity: 
\begin{equation}
  \label{eq:am4}
  P_{(+)}+P_{(-)}=I \ .
\end{equation}
(In particular, ~$P_{(+)}P_{(-)}=0$~ is the same as
(\ref{eq:am2}).)  
Thus, we  obtain the spectral decomposition: 
\begin{equation}
  \label{eq:am5}
  \hat{R}=(1-i)P_{(+)}+(1+i)P_{(-)} \quad = (1+i)I-2iP_{(+)} \quad
  = (1-i)I+2iP_{(-)}
\end{equation}

Note that although $\hat{R}$ is real, the roots in (\ref{eq:am2}) 
are complex and so are the projectors.

\subsection{Baxterisation}

We now apply the Baxterisation procedure for our case. 
First we introduce the following Ansatz (choosing a convenient
normalisation): 
\begin{equation}
  \label{eq:am6}
  \hat{R}(x)=I+c(x)\hat{R}
\end{equation}
and we try to find $c(x)$ such that $\hat{R}(x)$ would satisfy 
the parametrised Yang-Baxter equation:
\begin{equation}
  \label{eq:pybe}
\hat{R}_{(12)}(x)\hat{R}_{(23)}(xy)\hat{R}_{(12)}(y) -
\hat{R}_{(23)}(y)\hat{R}_{(12)}(xy)\hat{R}_{(23)}(x) = 0 
\end{equation}

With our Ansatz we actually have:
$$\hat{R}_{(12)}(x)\hat{R}_{(23)}(xy)\hat{R}_{(12)}(y) -
\hat{R}_{(23)}(y)\hat{R}_{(12)}(xy)\hat{R}_{(23)}(x) =$$
\begin{equation}
  \label{eq:am7}
  (c(x)+c(y)-c(xy))\bigl(\hat{R}_{(12)}-\hat{R}_{(23)}\bigr)
  +c(x)c(y)\bigl({\hat{R}^{2}}_{(12)}-{\hat{R}^{2}}_{(23)}\bigr) 
\end{equation}
Using (\ref{eq:am1}), i.e.,
$$\hat{R}^2=2(\hat{R}- I)$$ 
we obtain:
$$\hat{R}_{(12)}(x)\hat{R}_{(23)}(xy)\hat{R}_{(12)}(y) -
\hat{R}_{(23)}(y)\hat{R}_{(12)}(xy)\hat{R}_{(23)}(x)=$$
\begin{equation}
  \label{eq:am8}
  (c(x)+c(y)+2c(x)c(y)-c(xy))\bigl(\hat{R}_{(12)}-\hat{R}_{(23)}\bigr)
\end{equation}

Hence for Baxterisation one must have
\begin{equation}
  \label{eq:am9}
  c(x)+c(y)+2c(x)c(y)=c(xy)
\end{equation}
and then (\ref{eq:pybe}) holds. 

The solution of (\ref{eq:am9}) is:
\begin{equation}
  \label{eq:am10}
 2 c(x)= x^{p} -1 
\end{equation}

Is is interesting to note that (see, for example, \cite{isaev}, Sec.3.5),
the only change, 
as compared to $GL_{q}(N)$, is that one has a factor $2$ on the left
rather than 
$(q-q^{-1})$. Setting $p=-2$ and absorbing a free overall factor one
obtains, using also 
$$ 2I= \hat{R} +2\hat{R}^{-1}$$     
the elegant, symmetric form
\begin{equation}
  \label{eq:am11}
  \hat{R}(x)= (\sqrt{2}x)^{-1} \hat{R} +(\sqrt{2}x)\hat{R}^{-1} \;.
\end{equation}
Explicitly,
\begin{equation}
  \label{eq:baxtS03}
  \hat{R}(x) = \frac{1}{\sqrt{2x}}
  \left(
    \begin{array}{cccc}
      x+1 & 0 & 0 & 1-x \cr
      0 & x+1 & x-1 & 0 \cr
      0 & 1-x & x+1 & 0 \cr
      x-1 & 0 & 0 & x+1                                                   
    \end{array}
  \right)     \  .
\end{equation}

We shall also explore an Ansatz in terms of the projectors 
(instead of $\hat R$ as in (\ref{eq:am6})). 
Thus, for example, one may set 
\begin{equation}
  \label{eq:am12}
  \hat{R}(x)=I+a(x) P_{(+)}
\end{equation}
In this case the Baxterisation constraint turns out to be 
(compare with (\ref{eq:am9})):
\begin{equation}
  \label{eq:am13}
  a(xy)=\frac {a(x)+a(y)+a(x)a(y)}{1-\frac{1}{2}a(x)a(y)}
\end{equation}

The relation between $c(x)$ of (\ref{eq:am6}) and $a(x)$ of (\ref{eq:am12})
can be shown to be:
\eqna{eq:am18}
&&  a(x) =  \frac {2 c(x)} {i+(i-1)c(x)}  = 
  \frac {1+(1-i)c(x)}{1+(1+i)c(x)} -1 \\
&&  c(x) = \frac {i a(x)} {2 + (1-i) a(x)} 
\eena

Eq. (\ref{eq:am13}) is a special case of the functional equation:
\begin{equation}
  \label{eq:am14}
  a(xy)=\frac {a(x)+a(y)+a(x)a(y)}{1-k^{2}a(x)a(y)}
\end{equation}
which was  studied in a more general context in \cite{chakra3} 
where also  the solution of (\ref{eq:am14}) was found: 
\begin{equation}
  \label{eq:am15}
  a(x)=\frac {f(x)}{f(x^{-1})} -1
\end{equation}
where 
\begin{equation}
  \label{eq:am16}
  f(x)=x^{-1}-x \pm \sqrt {1-4k^2} (x+x^{-1})
\end{equation}
One obtains complex $a(x)$ due to the complex roots of (\ref{eq:am2}) and the
complex projectors in (\ref{eq:am3}).

In our case (cf. (\ref{eq:am14})) we have  
$$k^2=\frac{1}{2}$$
and hence
\begin{equation}
  \label{eq:am17}
  f(x)=x^{-1}-x \pm i (x+x^{-1})
\end{equation}
and it is sufficient to consider one, say the upper, sign - 
the lower sign will then correspond to the inversion  $x\to x^{-1}$. 
With this we obtain: 
\begin{equation}
  \label{eq:am17a}
  a(x)= \frac {x^2-1} {x^4+1} (1-x^2 + i(1+x^2) ) 
\end{equation}
Note that $a(x)$ is defined also for $x=0$, 
although the auxiliary function $f(x)$ is not.
Substituting this expression in (\ref{eq:am18}b) 
we recover (\ref{eq:am10}) for $p=-2$, i.e., the choice by 
which we obtained (\ref{eq:am11})) is not only 
a consistent one but is distinguished.  

On the other hand, if we replace $x$  by $x^{-p/2}$ 
in (\ref{eq:am17a}) and (\ref{eq:am18}b) 
then we shall recover exactly (\ref{eq:am10}). 
Such substitutions are possible since they do not change 
the character of equations (\ref{eq:am13}) and (\ref{eq:am14}).  

\subsection {Noncommutative Plane}

Here we shall find the noncommutative plane for our case. 
We shall apply the standard approach (cf. 
\cite{WZ}, also \cite{Schwenk}), however, not directly to 
the coordinates $(x_{1},x_{2})$ and the differentials
$(dx_{1},dx_{2})$ denoted as $\xi_{1},\xi_{2}$, but to the 
following complex linear combinations: 
\eqnn{eq:am19}
&&X_{1} = (x_{1}+ix_{2})\ , \quad X_{2} = (x_{1}-ix_{2}) \nn\\ 
&&  Z_{1}=(\xi_{1}+i\xi_{2})\ , \quad Z_{2}=(\xi_{1}-i\xi_{2})
\eea 
We introduce standard notation:
\begin{equation}
  \label{eq:am20}
  X\otimes X =
  \left(
 \begin{array}{c}
    X_{1}X_{1} \cr
    X_{1}X_{2} \cr
    X_{2}X_{1}\cr
    X_{2}X_{2}                                                  
   \end{array}
 \right)
\end{equation}
which we use also for $(X\otimes X) \mt (Z\otimes Z)$, 
$(X\otimes X) \mt (X\otimes Z)$, 
$(X\otimes X) \mt (Z\otimes X)$. 

According to the mentioned prescription for the 
consistent covariant calculus satisfying Leibniz rule, 
the commutation relations 
between the coordinates and differentials  are given as follows:
\eqna{ncpl} 
(\cP - I)\,   X\otimes X &=& 0 \\ 
(\cQ + I)\,   Z\otimes Z &=& 0 \\ 
\cQ\, (Z\otimes X) - X\otimes Z &=& 0 
\eena
where $\cP,\cQ$ are solutions of:
\eqn{qup} (\cP-I)\, (\cQ+I) = 0 
\end{equation}
The last equality in our situation 
means that ~$\cP-I,\cQ+I$~ are proportional 
to the projectors ~$P_+,P_-$~  or ~$P_-,P_+$, respectively. 
We choose:
\eqn{qupa} 
\cQ+I = \k\, P_+ \ , \qquad  \cP-I =  P_- 
\end{equation}
where ~$\k$~ is a complex proportionality constant 
and in the second relation we have set the proportionality constant 
equal to 1 due to the homogeneity of (\ref{ncpl}c).

Now the constraints on ~$x_{i}$~ and ~$\xi_{i}$~ can easily be
obtained. It turns out that 
the modular structure (\ref{ncpl}), expressed in terms of $x$ and
$\xi$, has real coefficients on choosing $\k$ imaginary, i.e.,
\begin{equation}
  \label{eq:am27}
  \k=ic \ , \quad c\in \RR \ .
\end{equation}
The definitions (\ref{eq:am19}) were introduced to assure this feature. The
final results are: 
\begin{equation}
  \label{eq:am28}
  {x_{1}}^2=x_{1}x_{2}, \quad {x_{2}}^2= -x_{2}x_{1}
\end{equation}
\begin{equation}
  \label{eq:am29}
  {{\xi}_{1}}^2=-{\xi}_{1}{\xi}_{2}, \quad {{\xi}_{2}}^2= {\xi}_{2}{\xi}_{1}
\end{equation}
and 
$$x_{1}{\xi}_{1}=(c-1){\xi}_{1}x_{1} + c{\xi}_{1}x_{2}$$
\begin{equation}
  \label{eq:am30}
  x_{1}{\xi}_{2}=(c-1){\xi}_{1}x_{2} + c{\xi}_{1}x_{1}
\end{equation}
$$x_{2}{\xi}_{1}=(c-1){\xi}_{2}x_{1} - c{\xi}_{2}x_{2}$$
$$x_{2}{\xi}_{2}=(c-1){\xi}_{2}x_{2} - c{\xi}_{2}x_{1}$$

For $c=1$ there is a supplementary simplification.

\section{$S14$}

\subsection{Spectral decomposition}

We take up now the second (of two) nonsingular non-triangular $R$-matrix 
in \cite{Hietarinta} which is not a deformation of the unit matrix:
\begin{equation}
\label{eq:S14}
 R_{S1,4}\ \equiv\ \left(
 \begin{array}{cccc}
 0 & 0 & 0 & q \cr
 0 & 0 & 1 & 0 \cr
 0 & 1 & 0 & 0 \cr
 q & 0 & 0 & 0 \cr 
 \end{array}
 \right)
\end{equation} 
where ~$q^2\neq 1$ (the case $q^2=1$ turned out \cite{ACDM2} 
to be equivalent to a special case of $GL_{p,q}(2)$). 
Here the braid matrix is 
\begin{equation}
  \label{eq:am31}
  \hat{R} = PR =
  \left(
 \begin{array}{cccc}
    0 &0 &0 &q \cr
    0 &1 &0 &0 \cr
    0 &0 &1 &0 \cr
    q &0 &0 &0                                                   
   \end{array}
 \right)
\end{equation}

The minimal polynomial identity $\hR$ satisfies is:
\begin{equation}
  \label{eq:am32}
  (\hat{R} -I) (\hat{R} -qI)(\hat{R} +qI) = 0  
\end{equation}

The projectors are: 
\eqnn{eq:am33}
&&P_{(0)}=\frac{1}{(1-q^2)}(\hat{R}-qI)(\hat{R}+qI) =
\left(
 \begin{array}{cccc}
  0 &0 &0 &0 \cr
  0 &1 &0 &0 \cr
  0 &0 &1 &0 \cr
  0 &0 &0 &0                                                   
 \end{array}
 \right) \nn\\
&&  P_{(+)}=\frac{1}{2q(q-1)}(\hat{R}+qI)(\hat{R}-I) =
  \frac{1}{2}
  \left(
 \begin{array}{cccc}
    1 &0 &0 &1 \cr
    0 &0 &0 &0 \cr
    0 &0 &0 &0 \cr
    1 &0 &0 &1                                                   
   \end{array}
 \right)
\\
&& P_{(-)}=\frac{1}{2q(q+1)}(\hat{R}-qI)(\hat{R}-I) =
\frac{1}{2} \left(
 \begin{array}{cccc}
  1 &0 &0 &-1 \cr
  0 &0 &0 &0 \cr
  0 &0 &0 &0 \cr
  -1 &0 &0 &1                                                   
 \end{array}
 \right)
\eea 
satisfying
$$P_{(i)}P_{(j)} = {{\delta}_{ij}}P_{(i)}$$
and
\begin{equation}
  \label{eq:am34}
  P_{(0)}+P_{(+)}+P_{(-)}=I
\end{equation}

Note that the projectors are independent of $q$. This leads to important
simplifications. The spectral decomposition is:
\begin{equation}
  \label{eq:am35}
  \hat{R}= \hR(q) = P_{(0)}+qP_{(+)}-qP_{(-)}=I+(q-1)P_{(+)}-(q+1)P_{(-)}
\end{equation}
from which it is obvious that:
\begin{equation}
  \label{eq:am36}
  \hat{R}^{-1}= \hR(q^{-1})
\end{equation}

\subsection{Baxterisation}

Here we introduce an Ansatz in terms of the projectors: 
\begin{equation}
  \label{eq:am37}
  \hat{R}(v,w)=I+vP_{(+)} + wP_{(-)}
\end{equation}
from which follows: 
\begin{equation}
  \label{eq:am38}
  \bigl(\hat{R}(v,w)\bigr)^{-1}=I-\frac{v}{1+v}P_{(+)} - \frac{w}{1+w}P_{(-)}
\end{equation}
The presence of two projectors above (as compared to one for
$S03$, cf. (\ref{eq:am12})) 
leads to a more elaborate structure analogous to that for $SO_{q}(N)$
(as compared to   $GL_{q}(N)$). 
We briefly display the analogy to the $SO_{q}(N)$ case
\cite{chakra2,chakra3} along with the 
drastic simplifications due to the special features of the projectors
noted above. 

Let us also introduce the following notations:
\eqn{notat} 
X_1 \equiv P_{(+)}\otimes I \ , \quad
X_2 \equiv I\otimes P_{(+)} \ , \quad 
Y_1 \equiv P_{(-)}\otimes I \ , \quad
Y_2 \equiv I\otimes P_{(-)} \ .\quad
\end{equation}
In which terms we have:
\eqn{notato} 
\hat{R}_{(12)}(v,w) = I\otimes I\otimes I + v X_1 + w Y_1 
 \ , \quad 
\hat{R}_{(23)}(v,w) = I\otimes I\otimes I + v X_2 + w Y_2 
\end{equation}

We obtain:
\eqnn{eq:am39}
&&\hat{R}_{(12)}(v,w)\hat{R}_{(23)}(v',w')\hat{R}_{(12)}(v'',w'') -
\hat{R}_{(23)}(v'',w'')\hat{R}_{(12)}(v',w')\hat{R}_{(23)}(v,w) =\nn\\
&&=\ (v+v''+vv''-v')S_{1}+(w+w''+ww''-w')S_{2}+  \nn\\
&&+\ vv'v''S_{5}+ww'w''S_{6}+
 \nn\\
&&+\ (vw'-v'w)J_{1}+(v''w'-v'w'')J_{2}+\nn\\
&&+\ vv'w''K_{1}+vw'v''K_{2}+wv'v''K_{3}+\nn\\
&&+\ ww'v''L_{1}+wv'w''L_{2}+vw'w''L_{3}
\eea
where 
\eqnn{eq:am40}
&&S_{1}= X_{1} - X_{2}\ ,\quad S_{2}= Y_{1} - Y_{2}  \ ,\quad
J_{1}= X_{1}Y_{2}-Y_{1}X_{2}\ , \quad J_{2}= Y_{2}X_{1}-X_{2}Y_{1}   \nn\\ 
&& S_{5}= X_{1}X_{2}X_{1}- X_{2}X_{1}X_{2}\ , \quad   S_{6}= Y_{1}Y_{2}Y_{1}-
Y_{2}Y_{1}Y_{2} \\   
&&K_{1}= X_{1}X_{2}Y_{1}-Y_{2}X_{1}X_{2}\ , \quad
K_{2}= X_{1}Y_{2}X_{1}-X_{2}Y_{1}X_{2}\ ,
\quad K_{3}=  Y_{1}X_{2}X_{1} - X_{2}X_{1}Y_{2}  \nn\\ 
&&  L_{1}= Y_{1}Y_{2}X_{1}-X_{2}Y_{1}Y_{2}\ , \quad
  L_{2}= Y_{1}X_{2}Y_{1}-Y_{2}X_{1}Y_{2}\ , 
  \quad L_{3}=  X_{1}Y_{2}Y_{1} - Y_{2}Y_{1}X_{2} \nn 
\eea 
First we note that from (\ref{eq:am35}) and (\ref{eq:am36}) it follows that
the right hand side of (\ref{eq:am39}) vanishes on setting 
\eqn{eq:am41}
v=v'=v''= q^{\pm1} -1, \qquad w=w'=w''= -(q^{\pm1} +1)
\end{equation}
However, not all these relations are necessary for the
vanishing of the RHS of (\ref{eq:am39}). More than this, the
successful Baxterisation would mean that this vanishing is
achieved with pairs $(v,w),(v',w'),(v'',w'')$ which are  
not identical (though satisfying some constraints in general). 

To accomplish this we first obtain a set of constraints relating
the members of (\ref{eq:am40}). For this we shall use the fact
that for any function $f(x)$, ${\hat R}^{\epsilon}$ denoting for
$\epsilon =\pm 1$ the matrices (\ref{eq:am35}) and
(\ref{eq:am36}) respectively, one has the well-known relations
\cite{FRT} (cf. also \cite{isaev}):
\eqnn{eq:am42}
&&f({\hat R}^{\epsilon}_{(12)}){\hat R}^{\epsilon}_{(23)}{\hat
  R}^{\epsilon}_{(12)}= 
{\hat R}^{\epsilon}_{(23)}{\hat R}^{\epsilon}_{(12)}f({\hat
  R}^{\epsilon}_{(23)}) \nn\\ 
&&  {\hat R}^{\epsilon}_{(12)}{\hat R}^{\epsilon}_{(23)}f({\hat
    R}^{\epsilon}_{(12)})= 
  f({\hat R}^{\epsilon}_{(23)}){\hat R}^{\epsilon}_{(12)}{\hat
    R}^{\epsilon}_{(23)} 
\eea 
For $f({\hat R}^{\epsilon}_{(12)})$ and $f({\hat
  R}^{\epsilon}_{(23)})$ one can 
choose $(X_1,Y_1)$ and $(X_2,Y_2)$ respectively and apply them in
(\ref{eq:am42}) successively. 

Exploiting systematically all the constraints implied by 
(\ref{eq:am42}) one obtains for the members of (\ref{eq:am40}): 
\eqnn{eq:am43}
&& S_{5}=\frac{1}{4} S_{1}, \qquad S_{6}=\frac{1}{4} S_{2} \\ 
&&K_{1}=-\frac{1}{4}(S_{2}+2J_{2}),\quad K_{3}=-\frac{1}{4}(
S_{2}+2J_{1}), \quad K_{2}=\frac{1}{4}(S_{2}+2(J_{1}+J_{2})) \nn\\
&&  L_{1}=-\frac{1}{4}(S_{1}-2J_{2}),\quad L_{3}=-\frac{1}{4}(
  S_{1}-2J_{1}), \quad L_{2}=\frac{1}{4}(S_{1}-2(J_{1}+J_{2})) \nn
\eea 

Hence the right hand side of (\ref{eq:am39}) becomes 
\begin{equation}
  \label{eq:am44}
  a_{1}S_{1}+a_{2}S_{2}+b_{1}J_{1}+b_{2}J_{2}
\end{equation}
where
\eqnn{eq:am45}
&&4a_{1}=4(v+v''+vv''-v')+ vv''v' -
w'wv''+v'ww''-w'w''v \nn\\ 
&&4a_{2}=4(w+w''+ww''-w')+ww''w' -
v'vw''+w'vv''-v'v''w \nn\\ 
&&2b_{1}= (w'v-v'w)(v''+w''+2) 
\nn\\ 
&&  2b_{2}=(w'v''-v'w'')(v+w+2) 
\eea 
The equations 
$$b_{1}=0,\qquad b_{2}=0$$
each has two factors that can be zero.  Setting the first factors
equal to zero  do not give satisfactory results when $a_{1}$ and
$a_{2}$ are equated to zero. So we set the second factors of both
$b_1$ and $b_2$ to zero:
\eqn{eq:am47}
v+w+2=0 \ , \qquad   v''+w''+2=0 \ ,
\end{equation}
which is consistent  with (\ref{eq:am41}). 
Taking into account (\ref{eq:am47}) leads to:
\eqn{eq:am48a}
2a_1 = 2a_2 = (v'+w'+2) (v+v''+vv'')
 \end{equation}
Consistently with (\ref{eq:am41}) the conditions $a_1 = a_2 =0$ 
lead to:
\begin{equation}
  \label{eq:am48}
  v'+w'+2  =0
\end{equation}

Thus, summarising, we have:
\begin{equation}
  \label{eq:am51}
  v+w=v'+w'=v''+w''=-2
\end{equation}
But apart from this the pairs ~$(v,w),(v',w'),(v'',w'')$~ are mutually 
{\it independent of each other}. 

Thus denoting
\begin{equation}
  \label{eq:am52}
  \hat{R}(q)=I+(q-1)P_{(+)}-(q+1)P_{(-)}
\end{equation}
satisfying
$$v+w =(q-1) -(q+1)=-2$$
one obtains
\begin{equation}
  \label{eq:am53}
  \hat{R}_{(12)}(q)\hat{R}_{(23)}(q')\hat{R}_{(12)}(q'') =
  \hat{R}_{(23)}(q'')\hat{R}_{(12)}(q')\hat{R}_{(23)}(q) 
\end{equation}
without further restrictions on the triplet$(q,q',q'')$.
Explicitly, $\hat R(q)$ is given by formula (\ref{eq:am31}).

One may choose
\begin{equation}
  \label{eq:am54}
  q=f(x), \qquad q''=f(y), \qquad q'=f(xy)
\end{equation}
to obtain a conventional Baxterisation. But more freedom is
implied in (\ref{eq:am53}). 

Let us consider some special cases. First we note that for 
$$w=w''=w'=0, \qquad a_{2}=b_{1}=b_{2}=0$$
and $a_{1}=0$ leads to 
\begin{equation}
  \label{eq:am55}
  v'=\frac{v+v''+vv''}{1-\frac{1}{4}vv''}
\end{equation}
Similarly for $v=v''=v'= 0$
\begin{equation}
  \label{eq:am56}
  w'=\frac{w+w''+ww''}{1-\frac{1}{4}ww''}
\end{equation}

Finally we note an amusing point. The general solution for the 
diagonaliser of $\hR$ contains arbitrary parameters. They can be so chosen
that the $\hR$ of $S03$ can be implemented to diagonalise the
$\hR$ of $S14$. Thus, setting (cf. (\ref{eq:am0}):
\begin{equation}
  \label{eq:am61}
  M =
  \frac{1}{\sqrt 2}
  \left(
 \begin{array}{cccc}
 1 &0 &0 &1 \cr
     0 &1 &-1 &0 \cr
    0 &1 &1 &0 \cr    
    -1 &0 &0 &1                                                   
   \end{array}
 \right)
\end{equation}
which fulfils:
$$ M\ ^tM ~=~ ^tM M = I $$ 
one obtains: 
\begin{equation}
  \label{eq:am62}
  M{\hat R}\ ^tM = diag (q,1,1,-q) 
\end{equation}
where $\hR$ is from (\ref{eq:am31}). (The mutually orthogonal rows of
$M$ may be permuted to reorder the diagonal elements of
(\ref{eq:am62}).) 

The situation is not reciprocal. Fixing suitably arbitrary
parameters one may construct a simple unitary matrix for
diagonalisation which we present for comparison. One can
use:
\begin{equation}
  \label{eq:am63}
  M' =
  \frac{1}{\sqrt 2}
  \left(
 \begin{array}{cccc}
 1 &0 &0 &i \cr
     0 &1 &-i &0 \cr
    0 &-i &1 &0 \cr    
    i &0 &0 &1                                                   
   \end{array}
 \right)
\end{equation}
which fulfils:
$$ M'\ {M'}^\dag ~=~ {M'}^\dag\ M' = I $$ 
to obtain: 
\begin{equation}
  \label{eq:am633}
  M'\ {\hat R'}\ {M'}^\dag = diag (1-i,1-i,1+i,1+i) 
\end{equation}
where ${\hat R'}$ is from (\ref{eq:am0}).

We mention only briefly here the utilities of explicit
construction of diagonalisers \cite{chakra3}, which here are $M$
and $M'$. They yield easily and directly the eigenvectors of
$\hR$, the latter being significant in related statistical
mechanical models. They also furnish new insights concerning
related noncommutative spaces.

\subsection{Noncommutative planes}

Using the notations and conventions used for $S03$, but using $x$ and
$\xi$ (not passing via $X$ and $Z$ as before), we use
(\ref{ncpl}) with $(X,Z)\to (x,\xi)$. As in (\ref{qupa}) we choose  
$$ \cP-I = P_- \Rightarrow P_{(-)}x\otimes x =0 $$
Following (\ref{qup}) we should have:
\begin{equation}
  \label{eq:am64}
  \bigl(\cQ+I \bigr) P_{(-)} =0
\end{equation}
from which follows that:
\begin{equation}
  \label{eq:am65}
\cQ = -I+2k_{+}P_{(+)}+k_{0} P_{(0)} 
\end{equation}
where $k_{+}, k_{0}$ are free parameters.
Thus, finally
  \eqna{eq:am66}
&&  P_{(-)}x\otimes x =0 \\
&&  ( 2k_{+}P_{(+)}+k_{0} P_{(0)}) 
{\xi}\otimes{\xi} =0\\
&&  x \otimes \xi=\bigl(-I+2k_{+}P_{(+)}+k_{0} P_{(0)}\bigr) {\xi} \otimes x 
\eena 
Note that from (\ref{eq:am66}b) follows:
$$ P_{(+)} {\xi}\otimes{\xi} =0 \ , \qquad 
P_{(0)} {\xi}\otimes{\xi} =0 $$
as $P_{(+)}$ and $P_{(0)}$ are orthogonal. 

The  relations resulting from (\ref{eq:am66}) are:
\eqn{ncx}
x_{1}^2 - x_{2}^2 =0
\end{equation}
\eqn{ncxi}
{\xi}_{1}^2 + {\xi}_{2}^2 =0 \ , \qquad 
{\xi}_{1}{\xi}_{2} = {\xi}_{2}{\xi}_{1} = 0 
\end{equation}
\eqnn{eq:am70}
&&x_{1}{\xi}_{1}= (k_{+} -1)\,{\xi}_{1}x_{1}+k_{+}\,{\xi}_{2}x_{2} \nn\\
&&x_{2}{\xi}_{2}= (k_{+} -1)\,{\xi}_{2}x_{2}+k_{+}\,{\xi}_{1}x_{1} \nn\\
&&x_{1}{\xi}_{2}= (k_{0} -1)\,{\xi}_{1}x_{2} \nn\\
&& x_{2}{\xi}_{1}= (k_{0} -1)\,{\xi}_{2}x_{1}
\eea 

Note that for $k_+=1$ and/or $k_0=1$ there are significant (even 
drastic for $k_0=1$) simplifications. 

\section{Conclusions}

In our first paper on exotic bialgebras \cite{ACDM1} we used 
triangular $R-$matrices which had multiple roots for $\hat R$ 
except for one case when it was no longer "exotic" but 
the non-standard Jordanian one. Various features of the 
Jordanian case have been studied (in a generalised biparametric form)
elsewhere \cite{chakra1}. For the remaining cases
mentioned above the multiple roots of $\hat R$ prevent
straightforward spectral decomposition. 

In the $S03,S14$ cases studied here  though there are no multiple
roots, one encounters in each case remarkable special features.

For $S03$ the projectors are complex. 
With the Ansatz for Baxterisation formulated directly in terms of
$\hat R$, (cf. (\ref{eq:am6})), a remarkable analogy with
the $GL_{q}$ case emerged. Following the crucial equations
(\ref{eq:am9}) and (\ref{eq:am10}), we noted that a coefficient
$(q-q^{-1})$ for the latter case is replaced by $2$ in ours.
Starting from an alternative Ansatz using a (complex) projector
it was shown again (cf.
(\ref{eq:am14}),(\ref{eq:am15}),(\ref{eq:am16})) how a particular 
case of a functional equation arising quite generally
\cite{chakra3} led to the solution.

For $S14$ one needs {\it
three} projectors for resolution of $I$ and the spectral
decomposition of $\hat R$. This introduces features analogous to
the more complicated $SO_{q}$ (and $Sp_{q}$) cases. One could
have arrived at our final results by working directly with the
explicit, numerical $8\times 8$ matrices for the tensored
projectors. Indeed, we have computed these matrices. But the
algebraic approach preferred here gives a deeper understanding of
the structure, emphasising both, the analogies and the
differences with the case of $SO_{q}(N)$.

\end{document}